\documentclass[12pt]{amsart}
\usepackage{amsfonts,amssymb}
\numberwithin{equation}{section}

 \def\numset#1{{\mathbb #1}}
 \def\setZ{\numset{Z}}
\def\setR{\numset{R}}
\def\setQ{\numset{Q}}

\makeatletter
\@namedef{subjclassname@2020}{%
  \textup{2020} Mathematics Subject Classification}
\makeatother

\theoremstyle{plain}
\newtheorem{Th}{Theorem}

 \theoremstyle{definition}

\newtheorem*{K}{Problem}

  \newcommand{\ab}[1]{\left\vert{#1}\right\vert}
 \newcommand{\zj}[1]{\left({#1}\right)}

\newcommand{\eq}[1]{\eqref{#1}}

\newcommand\ndiv{\nmid}

\newcommand\kong[2]{\equiv #1 \pmod{#2}}

\begin{document}

\baselineskip=17pt

\title{Beurling-integers with lacunarity}

\author{Imre Z. Ruzsa}
\address{Alfr\'ed R\'enyi Institute of Mathematics\\
     Budapest, Pf. 127\\
     H-1364 Hungary
}
\email{ruzsa@renyi.hu}
 \thanks{Supported by  NKFI grants K-129335, K-119528, KKP-133819.}


     \begin{abstract}
     We present examples of multiplicative semigroups of positive reals (Beurling's generalized integers)
with gaps bounded from below.
     \end{abstract}

  \maketitle

     \section{Introduction}

Let $G= \{g_1, g_2, \ldots \}$ be a sequence of real numbers, $1<g_1\leq g_2 \leq \ldots$ (generators) and
$B= \{ b_0, b_1, \ldots\}$, $b_0=1 <b_1 \leq b_2 \leq \ldots$ the sequence of products formed by elements of $G$.
If $G$ is the set of primes, $B$ will be the set of positive integers. The name honours Beurling,
who was the first to study analogs of the prime-number theorem for such systems.

 If $G$ is a set of
multiplicatively independent integers, $B$ will be a subset of positive integers, hence
$b_{i+1}-b_i \geq 1$. If furthermore $G$ contains all but finitely many primes, then $b_{i+1}-b_i$
will also be bounded from above. Lagarias \cite{lagarias97} 
proved that there is no other example consisting of integers,
and asked whether there is another example made of reals.

I conjecture that such a set does not exist. 

As a first step towards this conjecture, we show that a certain simple attempt to construct such a set
must fail, namely we cannot omit a small set of primes and replace them by non-integers.

\begin{Th}\label{nincs}
  Let $P$ be a set of primes such that
\begin{equation}\label{kimarad}
 %  E(x) =  | \{ p\leq x, p \notin P \}| = o(\sqrt{x}),
\sum_{p\notin P} 1/p < \infty \end{equation}
   and $\alpha\in\setR \setminus \setZ$, $\alpha>1$. With $G=P \cup \{ \alpha\}$ we have
    \[ \liminf b_{i+1}-b_i =0.\]
\end{Th}

On the other hand, we can add extra elements to a very thin set of primes.

\begin{Th}\label{kevesprim}
  Let $P$ be a set of primes such that
\[ 
    \sum_{p\in P} \frac{1}{\sqrt{p}} < \infty . \]
There exist numbers   $\alpha\in\setR \setminus \setZ$, $\alpha>1$ such that for $G=P \cup \{ \alpha\}$ we have
    $  b_{i+1}-b_i \geq  1$. The set of such numbers $\alpha$ has positive measure.
\end{Th}

The above theorem was stated as to form a contrast to Theorem \ref{nincs},
but in fact there is nothing special in the primes.

\begin{Th}\label{keves}
  Let $G'$ be a set of reals such that
  \begin{equation}\label{konvergens}
    \sum_{g\in G'} \frac{1}{\sqrt{g}} < \infty . \end{equation}
Let $B'$ be the sequence of products formed by elements of $G'$.
Assume that $b_{i+1}'-b_i' \geq \delta>0$ for all $i$.
There exist numbers   $\alpha\in\setR \setminus \setZ$, $\alpha>1$ such that for $G=G' \cup \{ \alpha\}$ we have
    $  b_{i+1}-b_i \geq \delta$. The set of such numbers $\alpha$ has positive measure.
\end{Th}

Unfortunately we cannot say much about sets of primes which are neither almost full
nor very thin. The metric appoach of Theorem \ref{keves} cannot be substantially improved.
We illustrate this by the example of squares, where conditon \eq{konvergens} ``just'' fails.

\begin{Th}\label{negyzet}
  Let $G'= \{ p^2\}$ be the set of prime squares, $B'= \{ n^2\}$ the set of squares.
 
There exist infinitely many numbers   $\alpha\in\setR \setminus \setZ$, $\alpha>1$ such that for $G=G' \cup \{ \alpha\}$ we have
    $  b_{i+1}-b_i \geq 1$. The set of such numbers $\alpha$ has measure 0.
\end{Th}

Call a set of Beurling-integers \emph{maximal lacunary}, if $ \inf  b_{i+1}-b_i>0$, but the inclusion
of any number to $G$ spoils this property.

\begin{K}
  How thin can a maximal lacunary Beurling-set be?  Is $B(x)= O(\sqrt{x})$ possible?
\end{K}

$x^{1/2+\varepsilon}$ is possible, as the following easy example shows.

\begin{Th}\label{kishijan}
  Let $1<c<2$, $G= \{ p^c\}$ be the set of $c$-th powers of primes, $B= \{ n^c\}$. This set is maximal lacunary.
\end{Th}

The densest example of a lacunary $B$ we could construct which is different from the integers is as follows.

\begin{Th}\label{pelda}
  There exists a set $G$ of irrational numbers such that
   \[ G(x) = \ab{\{ g\in G, g \leq x\}} > cx/\log x \]
and $b_{i+1}-b_i \geq 1$.
\end{Th}

\section{Proof of Theorem \ref{nincs}}

Let $E$ be the set of primes missing from $P$, and $R$ the set of integers composed only of
primes from $P$. We show that for every $\delta>0$ there are integers $x,y\in R$ such that
 \[ | \alpha^mx - y | < \delta .\]

 {\bf Case 1.} $\alpha$ is rational, say $\alpha=a/b$. We want to find  $x,y\in R$ with
  \[ \ab{ a^mx- b^my} < \delta b^m . \]

Fix $k$ so that $\delta b^m > 2$. Let $d=2$ if $ab$ is odd, $d=1$ otherwise. Fix odd numbers $u,v$ with
$a^mu-b^mv=d$. We will find $x,y$ in the form
 \[ x= u + 2zb^m, \ y = v+2za^m .\]
With such a choice we have $a^mx- b^my = d<  \delta b^m$. We need that $x,y$ be free of primes from $E$.
We shall estimate the number of integers $z<Z $ with this property.

For a prime $p\in E$, the divisibility $p| u + 2zb^m$ excludes at most one residue class 
modulo $p$. (Exactly one if $p \nmid b$ and none if $p|b$, since the assumption $a^mu-b^mv=d$ excludes
$p|(u,b)$.) For $p=2$ this divisibility cannot hold. Similarly  the divisibility $p|v+2za^m$ may
exclude a residue class, hence at least $p-2$ remain.

Write
 \[ \eta = \prod_{p\in E,  p>2} \zj{1-\frac{2}{p}} \]
 and select $T$ so that
  \[ \sum_{p\in E,  p>T} 1/p <  \eta/3 .\]
  Let $q= \prod_{p\in E, 2<p\leq T} p$. In each interval of length $q$ there are at least
   \[  \prod_{p\in E, 2<p\leq T} (p-2) \geq \eta q  \]
integers that avoid the excluded residue classes for every $p\leq T$. Up to $Z$
this is at least $\eta Z$ integers if $q|Z$.

Any prime divisor $p>T$ must be less than
 \[ \max ( u + 2zb^m,  v+2za^m) < cz, \]
 and excludes 2 residue classes which means at most $2(1+z/p)$ integers. There remain at least
  \[ \eta Z - \sum_{p\in E, T<p<cZ}   2(1+z/p) > \eta Z/3 - 2\pi(cZ) > 0 \]
for large $Z$.

{\bf Case 2.} $\alpha$ is irrational. Consider two convergents from the continued fraction development of $\alpha$ , say
 \[ \frac{a_k}{r_k} < \alpha <\frac{a_{k+1}}{r_{k+1}} . \]
 Any median
  \[ \mu =  \frac{xa_k+ya_{k+1}}{xr_k+yr_{k+1}}, \ x,y>0   \]
satisfies
\[ \frac{a_k}{r_k} < \mu <\frac{a_{k+1}}{r_{k+1}} , \]
hence
 \[ |\alpha-\mu| < \frac{a_{k+1}}{r_{k+1}} - \frac{a_k}{r_k} = \frac{1}{r_kr_{k+1}} . \]

 We try to find $x,y$ so that the numerator and denominator of $\mu$  be free of primes from $E$ and
\begin{equation}\label{kozelvan}
\ab{ \alpha ( {xr_k+yr_{k+1}}) - ({xa_k+ya_{k+1}})} < \frac{xr_k+yr_{k+1}} {r_kr_{k+1}} < \delta. \end{equation}
For the last inequality to hold we require
\begin{equation}\label{korlatok}
x<X= \delta r_{k+1}/2, \ y <Y = \delta r_k/2. \end{equation}

First we fix the  parity of $x$ and $y$ to make  the numerator and denominator odd.
If $a_kr_k$ is odd, we set $2|y$,  $2\nmid x$. If $a_{k+1}r_{k+1}$ is odd, we set $2\nmid y$,
 $2|x$. If neither happens, then we set $2\nmid y$ and  $2\nmid x$. The fact
that $a_{k+1}r_k - a_kr_{k+1} =1 $ ensures that this works. 

Given $y$, for a prime $p>2$ the divisibility $p|{xa_k+ya_{k+1}}$ means a single residue class modulo $p$ if 
$p \nmid a_k$. It is impossible if $p|a_k$ and $p \ndiv y$, and it always holds if $p|(y, a_k)$.
 Similarly, the divisibility $p|{xr_k+yr_{k+1}}$ 
means a single residue class modulo $p$ if 
$p \nmid r_k$, 
it is impossible if $p|r_k$ 
and $p \ndiv y$, and it always holds if $p|(y, r_k)$.
 That is, at most two residue classes are excluded
modulo $p$ unless $p|(y, a_k)$ or $p|(y, r_k)$. As we have little control over prime divisors of
$a_k$ and $r_k$, we will require that $y$ be free of prime divisors from $E$ up to a limit.

Write
 \[ \eta = \frac{1}{2}\prod_{p\in E, p>2} \zj{1-\frac{2}{p}}, \ \eta' = \frac{1}{2}\prod_{p\in E, p>2} \zj{1-\frac{1}{p}} \]
 and select $T$ so that
  \[ \sum_{p\in E, p>T} 1/p <  \delta\eta\eta'/10 .\]
  Let $q=2 \prod_{p\in E, 2<p\leq T} p$. 

In each interval of length $q$ there are at least
   \[  \prod_{p\in E, 2<p\leq T} (p-1) \geq \eta' q  \]
integersfree of prime divisors  $p\in E$, $p\leq T$.
 Up to $Y$
this is at least $\eta'Y-q$ integers. For each such $y$,
in each interval of length $q$ there are at least
   \[  \prod_{p\in E, 2<p\leq T} (p-2) \geq \eta q  \]
integers that avoid the excluded residue classes for every $p\leq T$. Up to $X$
this is at least $\eta X-q$ integers. This leaves us with at least
 \[ (\eta X-q)(\eta'y-q) > \delta^2\eta\eta' r_kr_{k+1}/5 \]
possible pairs $(x,y)$.

Consider prime divisors $p>T$. The integers, which should not be divisible by these primes,
are all less than
 \[ Xa_k + Ya_{k+1} < (\delta/2)(a_k r_{k+1} +r_k a_{k+1})< U=2\delta r_k  r_{k+1} ; \]
hence this is also a bound for $p$. The numbers $ xa_k+ya_{k+1}$ are all distinct by the coprimality
of $a_k$ and $a_{k+1}$, and so are the numbes $xr_k+yr_{k+1}$, but we cannot exclude that the two kinds
overlap. Hence an upper estimate for the number of pairs $x,y$ with an illegal divisibility
is $2(U/p+1)$. Summing this for all $p<U$ we obtain
 \[ \sum_{T<p<U}2(U/p+1) <2U \sum_{p\in E, p>T} 1/p + 2\pi(U) <  \delta^2\eta\eta' r_kr_{k+1}/5   \]
if $r_k$ is large enough.

 \section{Proof of Theorems \ref{keves} and \ref{kevesprim}} 

We need to find numbers $\alpha$ such that $\ab { \alpha^km - \alpha^jn} \geq \delta$ for all
$m,n\in G'$ and positive integers $j<k$. Since for $j \leq k$ we have
\[ \ab { \alpha^km - \alpha^jn} = \alpha^j  \ab { \alpha^{k-j}m - n} \geq\ab { \alpha^{k-j}m - n} , \]
it is sufficient to consider the case $j=0$. 

 We will show that the measure of such $ \alpha$ in the interval 
 $[e^t, e^{2t}]$ is positive for sufficiently large $t$.

 The event we want to avoid is  $\ab { \alpha^km - n} < \delta$, which can be rewritten as
  \[ \alpha^k \frac{m}{n} \in \zj{1-\frac{\delta}{n},1+\frac{\delta}{n }} .  \]
Note that  $\ab { \alpha^km - n} < \delta$ implies $n>\alpha m-\delta$, whence $n>2\delta$ and $n>\alpha m/2>m$, assuming
that $\alpha>3\delta$ which holds for $t> \log 3\delta$.

   We take logarithms to infer, with the notation $\beta= \log \alpha$, that
    \[ k\beta + \log m - \log n \in (-2\delta/n, \delta/n) ,\]
    that is,
      \[ \beta \in \frac{\log n - \log m}{k} + \zj{ \frac{-2\delta}{kn},  \frac{\delta}{kn}} . \]
  To estimate the measure of such numbers $\beta$ we add $ 3\delta/(kn)$ for all triplets
  $m,n,k$ such that the above interval intersects the interval $[t, 2t]$. If  $t>4\delta$,
 this intersection implies
   \[ \frac{\log n - \log m}{k} \in (t/2, 3t) .\]
   Hence
    \[  \frac{\log n - \log m}{3t} <k < 2 \frac{\log n - \log m}{t}.  \]
The ratio of the upper and lower bounds is 6, hence the sum of $1/k$ in this interval is less thn
$c=1+ \log 6$. Consequently the sum of  $ 3\delta/(kn)$ for all triplets
$m,n,k$ is at most the sum of $3c\delta/n$ for all possible pairs $m,n$. These pairs staisfy $n>\alpha m/2> e^tm/2$, so
 \[ \sum_{m,n} \frac{1}{n} < 2e^{-t/2} \sum_{m,n\in B'} \frac{1}{\sqrt{mn}} = 2e^{-t/2} \zj{ \sum_{m\in B' } \frac{1}{\sqrt{m}}}^2 . \]
 This series is convergent, indeed
  \[  \sum_{m\in B'} \frac{1}{\sqrt{m}} = \prod_{g\in G'} \zj{1+ \frac{1}{\sqrt{g}-1}} < \infty \]
by assumption \eq{konvergens}.

The estimate we found for the measure of bad $\beta$ is
 \[ 6c\delta e^{-t/2} \zj{ \sum_{m\in B' } \frac{1}{\sqrt{m}}}^2 ,\]
which is less than $t$, the length of the interval for large enough $t$.

\section{Proof of Theorem \ref{negyzet}}

Let $q$ be a squarefree integer, $a,b$ positive integers and
 \[ \alpha = \bigl(a\sqrt{q} +b\bigr)^2 .\]
We show that for these numbers $B$ has the lacunarity property.

 The elements of $B$ are numbers of the form $\alpha^km^2$, and we need to show that
  \[ \ab { \alpha^km^2 - \alpha^jn^2} \geq 1 . \]
Since for $j \leq k$ we have
\[ \ab { \alpha^km - \alpha^jn} = \alpha^j  \ab { \alpha^{k-j}m - n} \geq\ab { \alpha^{k-j}m - n} , \]
it is sufficient to consider the case $j=0$.

Put $\beta =  \bigl(a\sqrt{q} +b\bigr)^k$. This number is of the form
 \[ \beta = u\sqrt{q} +v \]
 with positive integers $u,v$. Now we have
  \[ \alpha^km^2-n^2 = (\beta m)^2-n^2 = \bigl( vm+n+um\sqrt{q}\bigr)\bigl( vm-n+um\sqrt{q}\bigr) \] \[
= \frac{ \bigl( vm+n+um\sqrt{q}\bigr)}{ \bigl( vm-n-um\sqrt{q}\bigr)} \zj{ (vm-n)^2 -(um)^2q} .
 \]
The enumerator exceeds the absolute value of the denominator, and the second factor is a nonzero
integer, so the absolute value of the expression is $>1$.

 Now we show that for such numbers $\sqrt{\alpha}$ are badly approximable. Assume that it is well approximable, that is,
  for every $\varepsilon>0$ there are integers
 $a,b$ such that
  \[ \ab{\sqrt{\alpha}-\frac{a}{b}} < \frac{\varepsilon}{b^2} . \]
  Clearly $a<2\sqrt{\alpha}b$ and then
    \[ \ab{\alpha b^2-a^2} = (\sqrt{\alpha}b-a)(\sqrt{\alpha}b+a) < 3\varepsilon\sqrt{\alpha}.\]
  Badly approximable numbers have measure 0 by a theorem of  Hinchin \cite{hincsin}. 

\section{Proof of Theorem \ref{kishijan}}

Try to include a number $\alpha$. Take integers $a,b$ such that
 \[ \ab{\alpha^{1/c}- \frac{a}{b}} < \frac{1}{b^2} . \]
 From the mean value theorem we see that
  \[ \frac{\alpha b^c-a^c}{\alpha^{1/c}b-a} = cz^{c-1} \]
  with some $z$ between $\alpha^{1/c}b$ and $a$, so $z=O(b)$. Hence
   \[ \alpha b^c-a^c = O(b^{c-2}) \]
can be arbitrarily small. 

\section{Proof of Theorem \ref{pelda}}

We give two examples, one with quadratic irrationals and the other with transcendental numbers.
Both arise from a subset of primes through a transformation.

{\bf Example 1: quadratic.}

Take those odd primes that split in $\setQ[\sqrt{2}]$. They are the primes $p \kong{\pm 1}{8}$ (about half of 
the primes). For such a prime there are positive integers $x,y$ such that
 \[ \pm p = x^2-2y^2 = (x- y\sqrt{2})(x+ y\sqrt{2}).  \]
Put $f(p) = \min (x+ y\sqrt{2})$ over all such representations. This satisfies 
$f(p) < C\sqrt{p}$ with some constant $C$. This can be seen by comparing the minimal representation
with the one obtained by $x'= |x-2y|$, $y'=|y-x|$ which corresponds to a multiplication 
by the unit $1-\sqrt{2}$ of $\setQ[\sqrt{2}]$. (It is not difficult to calculate the best value of $C$,
but not too important for this argument.) Extend $f$ multiplicatively to all integers composed
exclusively of primes $p \kong{\pm 1}{8}$. For every such integer $n$ we have
 \[ f(n) =x+ y\sqrt{2}, \ x,y>0, \  x^2-2y^2 =n . \]

Put $g(n) = f(n)^2$. Our generators will be the numbers $g(p)$ for $p$ prime, $B$ will be the set
of values of $g(n)$ for the above described special $n$. As $g(p)<C^2p$ and half of the primes are used,
$G(x) > cx/\log x$ holds for large $x$ with $c=1/(2C^2)$. 

Now we show that $| g(m)-g(n)| > 1$ for $m\neq n$. Let
 \[ f(m) = u+v\sqrt{2}, \ f(n) =x+ y\sqrt{2} . \]
 We have
  \[ f(m)^2 - f(n)^2  = \zj{ (u+x) + (v+y)\sqrt{2} } \zj{ (u-x) + (v-y)\sqrt{2}} \] 
\[ =  \frac{ (u+x) + (v+y)\sqrt{2}}{ (u-x) - (v-y)\sqrt{2}} \zj{ (u-x)^2 -2 (v-y)^2}
 . \]

The enumerator exceeds the absolute value of the denominator, and the second factor is a nonzero
integer, so the absolute value of the expression is $>1$.

The similarity to the proof of Theorem \ref{negyzet} hints that the two arguments could be combined,
and the above example can be extended by including squares of integers. However, this does not
 substantially increase the size of $B(x)$ and $G(x)$.

{\bf Example 2: transcendental.}

Consider primes $p \kong{1}{4}$. Write $p=a^2+b^2$ with $0<a<b$ and let
 \[ \rho(p) = ia+b = \sqrt{p} e^{ih(p)}, \ 0 < h(p) < \pi/2 .\]
Here $\rho(p)$ is one of the Gaussian primes in the decomposition of $p$ inthe ring of Gaussian integers.
Extend $\rho$ multiplicatively to the product of such primes, that is, odd integers that can be written
as a sum of two squares.

Since together with a Gaussian prime its conjugate is never selected, the numbers $\rho(n)$ for $n \neq 1$,
and $\rho(m)/\rho(n)$ for $m \neq n$ will never be real. Indeed, $\rho(m)/\rho(n)$ is a product of our selected primes
with (positive and negative) exponents, and its conjugate can be obtained by taking the conjugate primes,
and by the unicity of prime factorization these are different numbers.

Given a prime $p$ let $f(p)=h(p)+2k\pi$ with the integer $k$ chosen so that
$\log p < f(p) < \log p + 2\pi$. Extend $f$ additively. We will always have
 \[ e^{if(n)} = \frac{ \rho(n)}{\sqrt{n}} . \]
   Finally we put
    \[ g(p) = e^{f(p)} < e^{2\pi} p . \]
    These numbers form the set $G$, and (since again half of the primes was used)
     \[ G(x) > cx/\log x, \ c=e^{-2\pi}/2 \]
for large $x$. $B$ is the set of values of the multiplicative extension of $g$.
Since $g(n)$ is one of the values of $ \zj{ \rho(n)/\sqrt{n}}^{-i}$, it is transcendental
by the Gelfond-Schneider theorem, see for insance \cite{burgertubbs}.

We show the lacunarity property. For $m \neq n$ consider the triangle in the integer lattice
with vertices $0, \rho(m), \rho(n)$. Since it is a nondegenerate triangle, its area is at least 1/2,
on the other hand it is exactly
 \[  \frac{1}{2}  \sqrt{mn} \ab{ \sin \zj{f(m)-f(n)}} . \]
We infer that
\[  \ab{ \sin \zj{f(m)-f(n)}} \geq \frac{1}{\sqrt{mn}} . \]
Finally
 \[ g(m)-g(n) = e^{f(m)} - e^{f(n)} = e^{\frac{f(m)+f(n)}{2}} \zj{  e^{\frac{f(m)-f(n)}{2}} -  e^{\frac{f(n)-f(m)}{2}} } .\]
 The first factor is
  \[ \sqrt{ g(m) g(n)} \geq \sqrt{mn} . \]
  To estimate the second note that
   \[ \ab {e^x - e^{-x}} > 2 |x|  > | \sin (2x) |, \]
so it exceeds $  \ab{ \sin \zj{f(m)-f(n)}} $ which was shown to exceed $1/\sqrt{mn}$.

% \section{Remarks}
\bf{Acknowledgement. } This work was inspired by converstions with Szil\'ard G. R\'ev\'esz.

\providecommand{\bysame}{\leavevmode\hbox to3em{\hrulefill}\thinspace}
\providecommand{\MR}{\relax\ifhmode\unskip\space\fi MR }
\providecommand{\MRhref}[2]{%
  \href{http://www.ams.org/mathscinet-getitem?mr=#1}{#2}
}
\providecommand{\href}[2]{#2}

 %   \bibliographystyle{amsplain}
%        \bibliography{cimek}
%        \begin{thebibliography}{9}

  %     \end{thebibliography}
% \end{thebibliography}

 \end{document}